\documentclass[10pt]{article}
\usepackage{amsmath,amssymb,amsfonts}
\newtheorem{theorem}{Theorem}
\newtheorem{lemma}{Lemma}
\newtheorem{corollary}{Corollary}

\date{}

\title{On the Largest Part Size and Its Multiplicity of a Random Integer Partition}
\bigskip

\author{{\bf Ljuben Mutafchiev}\\
American University in Bulgaria, 2700 Blagoevgrad, Bulgaria \\ and
Institute of Mathematics and Informatics of the \\ Bulgarian
Academy of Sciences
\\ \tt {ljuben@aubg.edu}}

\begin{document}
\maketitle

\begin{abstract}
Let $\lambda$ be a partition of the positive integer $n$ chosen
uniformly at random among all such partitions. Let
$L_n=L_n(\lambda)$ and $M_n=M_n(\lambda)$ be the largest part size
and its multiplicity, respectively. For large $n$, we focus on a
comparison between the partition statistics $L_n$ and $L_n M_n$.
In terms of convergence in distribution, we show that they behave
in the same way. However, it turns out that the expectation of $L_n
M_n -L_n$ grows as fast as $\frac{1}{2}\log{n}$. We obtain a
precise asymptotic expansion for this expectation and conclude
with an open problem arising from this study.
\end{abstract}

\vspace{.5 cm}



 \vspace{.2cm}

\section{Introduction and statement of the main results }

Partitioning integers into summands (parts) is a subject of
intensive research in combinatorics, number theory and statistical
physics. If $n$ is a positive integer, then by a partition
$\lambda$ of $n$, we mean the
representation
\begin{equation}\label{lambd}
\lambda: \quad n=\lambda_1+\lambda_2+.+\lambda_k, \quad k\ge 1,
\end{equation}
where the positive integers $\lambda_j$ satisfy
$\lambda_1\ge\lambda_2\ge...\ge\lambda_k$. Let $\Lambda(n)$ be the
set of all partitions of $n$ and let $p(n)=\mid\Lambda(n)\mid$ (by
definition $p(0)=1$ regarding that the one partition of $0$ is the
empty partition). The number $p(n)$ is determined asymptotically
by the famous partition formula of Hardy and Ramanujan \cite{HR18}:
$$
p(n)\sim\frac{1}{4n\sqrt{3}}\exp{\left(\pi\sqrt{\frac{2n}{3}}\right)},
\quad n\to\infty.
$$
A precise asymptotic expansion for $p(n)$ was found by Rademacher
\cite{R37} (see also \cite[Chapter~5]{A76}). Further on, for fixed integer
$n\ge 1$, a partition $\lambda\in\Lambda(n)$ is selected uniformly
at random. In other words, we assign the probability $1/p(n)$ to
each $\lambda\in\Lambda(n)$. We denote this probability measure on
$\Lambda(n)$ by $\mathbb{P}$. Let $\mathbb{E}$ be the expectation
with respect to $\mathbb{P}$. In this way, each numerical
characteristic of $\lambda\in\Lambda(n)$ can be regarded as a
random variable, or, a statistic produced by the random generation
of partitions of $n$.

In this paper, we focus on two statistics of random integer
partitions $\lambda\in\Lambda(n)$: $L_n=L_n(\lambda)=\lambda_1$,
which is the largest part size in representation (\ref{lambd}) and
$M_n=M_n(\lambda)$, equal to the multiplicity of the largest part
$\lambda_1$ (i.e., $M_n(\lambda)=m,1\le m\le k-1$, if
$\lambda_1=...=\lambda_m>\lambda_{m+1}\ge...\ge\lambda_k$ in
(\ref{lambd}), and $M_n(\lambda)=k$, if
$\lambda_1=...=\lambda_k$).

Each partition $\lambda\in\Lambda(n)$ has a unique graphical
representation called Ferrers diagram \cite[Chapter~1]{A76}. It
illustrates (\ref{lambd}) by the two-dimensional array of dots,
composed by $\lambda_1$ dots in the first (most left) row,
$\lambda_2$ dots in the second row,..., and so on, $\lambda_k$
dots in the $k$th row. Therefore, a Ferrers diagram may be
considered as a union of disjoint blocks (rectangles) of dots
whose side lengths represent the part sizes and their multiplicities
of the partition $\lambda$, respectively. For instance Figure 1 illustrates the
partition $7+5+5+5+4+2+1+1+1$ of $31$ in which $L_n=7$ and
$M_n=1$.
$$
\begin{array}{ccccccc}
\bullet & \bullet & \bullet & \bullet & \bullet & \bullet &
\bullet \\
\bullet & \bullet & \bullet & \bullet & \bullet \\
\bullet & \bullet & \bullet & \bullet & \bullet \\
\bullet & \bullet & \bullet & \bullet & \bullet \\
\bullet & \bullet & \bullet & \bullet \\
\bullet & \bullet \\
\bullet \\
\bullet \\
\bullet \\
\end{array}
$$

$$
\quad {\bf Figure \quad 1}
$$

The earliest asymptotic results on random integer partition
statistics has been obtained long ago by Husimi \cite{H38} and
Erd\"{o}s and Lehner \cite{EL41}. Husimi has derived an asymptotic
expansion for $\mathbb{E}(L_n)$ in the context of a statistical
physics model of a Bose gas. Erd\"{o}s and Lehner were apparently
the first who have studied random partition statistics in terms of
probabilistic limit theorems. In fact, they showed that
\begin{equation}\label{el}
\lim_{n\to\infty}\mathbb{P}\left(\frac{L_n}{\sqrt{n}}-\frac{1}{2c}
\log{n}\le u\right) =H(u),
\end{equation}
where
\begin{equation}\label{gumbel}
H(u)=\exp{\left(-\frac{1}{c} e^{-cu}\right)}, \quad
-\infty<u<\infty
\end{equation}
and
\begin{equation}\label{c}
c=\frac{\pi}{\sqrt{6}}.
\end{equation}
Husimi's asymptotic result was subsequently reconfirmed by Kessler
and Livingston \cite{KL76}. Higher moments of $L_n$ were studied in
\cite{R74}. A general method providing asymptotic expansions of
expectations of integer partition statistics was recently proposed
by Grabner et al. \cite{GKW14}. Among the numerous examples, they derived a
rather complete asymptotic expansion for $\mathbb{E}(L_n)$, namely,
\begin{eqnarray}\label{eel}
& & \mathbb{E}(L_n) =\frac{\sqrt{n}}{2c}(\log{n}+2\gamma-2\log{c})
+\frac{\log{n}}{2c^2} +\frac{1}{4} \nonumber \\
& & +\frac{1+2\gamma-2\log{c}}{4c^2}
+O\left(\frac{\log{n}}{n}\right), \quad n\to\infty,
\end{eqnarray}
where $c$ is given by (\ref{c}) and $\gamma=0.5772...$ denotes the
Euler's constant (see \cite[Proposition~4.2]{GKW14}). Notice that by
conjunction of the Ferrers diagram the largest part and the total
number of parts in a random partition of $n$ are identically
distributed for any $n$. The sequence
$\{p(n)\mathbb{E}(L_n)\}_{n\ge 1}$ is given in \cite{S} as A006128.

There are serious reasons to believe that the multiplicity $M_n$ of the largest part of a random
partition of $n$ behaves asymptotically in a much simpler way than many other
partition statistics.  Grabner and Knopfmacher \cite{GK06} used Erd\"{o}s-Lehner limit theorem
(\ref{el}) to establish that
 \begin{equation}\label{em}
 \lim_{n\to\infty}\mathbb{E}(M_n)=1.
\end{equation}
In addition,
among many other important asymptotic
results, Fristedt, in his remarkable paper \cite{F93}, showed that, with
probability tending to $1$ as $n\to\infty$, the first $m_n$
largest parts in a random partition of $n$ are distinct if
$m_n=o(n^{1/4})$. Hence it may not be that $L_n$ constitutes the main
 contribution to $n$ by a single part size and some smaller part
 sizes may occur with sufficient multiplicity so that the products
 of these part sizes with their multiplicities could be much
 larger than $L_n$. In terms of the Ferrers diagram this means
 that its first block has typically smaller area than several next
 block areas with larger heights (multiplicities of parts).

Our aim in this paper is to study the asymptotic behavior of the
area $L_n M_n$ of the first block in the Ferrers diagram of a random
partition of $n$. We show some similarities and
differences between the single part size $L_n$ and its
corresponding block area $L_n M_n$. As a first step, we obtain a
distributional result for $M_n$ that confirms the limit in
(\ref{em}).

\begin{theorem} For any $n\ge 1$, we have
\begin{equation}\label{mnexact}
\mathbb{P}(M_n=1)=\frac{p(n-1)}{p(n)}.
\end{equation}
In addition, if $n\to\infty$, then
\begin{equation}\label{mnlimit}
\mathbb{P}(M_n=1) =1-\frac{c}{\sqrt{n}}+\frac{1+c^2/2}{n}+
O(n^{-3/2}),
\end{equation}
where the constant $c$ is given by (\ref{c}).
\end{theorem}

Combining the Erd\"{o}s-Lehner limit theorem (\ref{el}) with
(\ref{mnlimit}), one can easily observe that the limiting distributions
 of $L_n$ and $L_n M_n$ coincide under the same
normalization.

\begin{corollary} For any real $u$, we have
$$
\lim_{n\to\infty}\mathbb{P}\left(\frac{L_n M_n}{\sqrt{n}}
-\frac{1}{2c}\log{n}\le u\right) =H(u),
$$
where $H(u)$ and $c$ are given by (\ref{gumbel}) and (\ref{c}),
respectively.
\end{corollary}

Although $L_n$ and $L_n M_n$ follow the same limiting
distribution, the difference in means $\mathbb{E}(L_n
M_n)-\mathbb{E}(L_n)$ grows as fast as $\frac{1}{2}\log{n}$. A
complete estimate is given in the following

\begin{theorem} If $n\to\infty$, then, as $n\to\infty$,
$$
\mathbb{E}(L_n M_n) =\mathbb{E}(L_n) +\frac{1}{2}\log{n}-C
+O(1/\log{n}),
$$
where $C=\log{c}+1-\gamma=0.67165...$ and $\mathbb{E}(L_n)$ and
$c$ are given by (\ref{eel}) and (\ref{c}), respectively.
\end{theorem}

{\it Remark 1.} The sequence $\{p(n)\mathbb{E}(L_n M_n)\}_{n\ge
1}$ is given as A092321 in \cite{S}.

The proofs of Theorems 1 and 2 are based on generating function
identities established in \cite{ABBKM16} that involve products of the
form $P(x)G(x)$, where $P(x)$ is the Euler partition generating
function
\begin{equation}\label{euler}
P(x):=\sum_{n=0}^\infty p(n)x^n=\prod_{j=1}^\infty (1-x^j)^{-1}
\end{equation}
and $G(x)$ is a function which is analytic in the open unit disk
and does not grow too fast as $x\to 1$. Our  asymptotic expansions
in (\ref{mnlimit}) and Theorem 2 are obtained using a general
asymptotic result of Grabner et al. \cite{GKW14} for the $n$th
coefficient $x^n[P(x)G(x)]$. In the proof of Theorem 2 we also
apply a classical way of estimating the growth of a power series
around its main singularity.

We organize the paper as follows. Section 2 contains some
auxiliary facts related to generating functions and the asymptotic
analysis of their coefficients. In Section 3 we present the  proofs
of the main results. Finally. in Section 4 we conclude with an
open problem on the position of $L_n M_n$ in the sequence of
ordered block areas of a random Ferrers diagram.

\section{Preliminaries: Generating Functions and an Asymptotic
Scheme}

We start with two generating function identities. In the next two
lemmas $P(x)$ is the generating function given by (\ref{euler})
and by definition $\prod_1^0:=1$.

\begin{lemma} (i) For any positive integer $m$, we have
$$
\sum_{n=1}^\infty p(n)\mathbb{P}(M_n=m)x^n =x^m\prod_{j\ge
m}(1-x^j)^{-1} =P(x)x^m\prod_{j=1}^{m-1} (1-x^j).
$$

(ii) We have
$$
\sum_{n=1}^\infty p(n)\mathbb{E}(L_n M_n)x^n =\sum_{k=1}^\infty
\frac{kx^k}{1-x^k}\prod_{j=1}^k (1-x^j)^{-1} =P(x)F(x),
$$
where
\begin{equation}\label{efx}
F(x)=\sum_{k=1}^\infty\frac{kx^k}{1-x^k}\prod_{j=k+1}^\infty
(1-x^j).
\end{equation}
\end{lemma}

{\it Sketch of the proof.} Part (i) is the last conclusion of
Theorem 2.3 from \cite{ABBKM16}. Part (ii) is given in A092321 of \cite{S}.
It also follows from Proposition 4.1 in \cite{ABBKM16}. $\Box$

We shall essentially use the main result from \cite[Theorem~2.3]{GKW14}.
We present here only slight modifications of those parts of this theorem that we will need for
our further asymptotic analysis. Furthermore, by $\log{x}$ we
denote the main branch of the logarithmic function that satisfies
the inequality $\log{x}<0$ if $0<x<1$.

\begin{lemma} Suppose that, for some constants $K>0$ and $\eta<1$,
the function $G(x)$ satisfies
\begin{equation}\label{fast}
G(x)=O(e^{K/(1-\mid x\mid)^\eta}), \quad \mid x\mid\to 1.
\end{equation}

(i) Let $G(e^{-t})=at^b+O(\mid f(t)\mid)$ as $t\to 0,$
$\Re{t}>0$, where $b\ge 0$ is an integer and $a$ is real number.
Then, we have
\begin{eqnarray}
& & \frac{1}{p(n)} x^n[P(x)G(x)]
=a\left(\frac{2\pi}{\sqrt{24n-1}}\right)^b \frac{s}{s-1}
\sum_{j=0}^{b+1}\frac{(b+j+1)!}{j!(b+j-1)!}\left(-\frac{1}{2s}\right)^j
\nonumber \\
& & +O(e^{-2s}) +O\left(e^{-n^{1/2-\epsilon}}
+f(c/\sqrt{n}+O(n^{-1/2-\epsilon})\right) \nonumber
\end{eqnarray}
for any $\epsilon\in (0,(1-\eta)/2)$, where
\begin{equation}\label{es}
s=\sqrt{\frac{2\pi^2}{s}\left(n-\frac{1}{24}\right)}
=2c\sqrt{n-\frac{1}{24}}
\end{equation}
and $c$ is given by (\ref{c}).

(ii) Suppose that $G(x)$ satisfies condition (\ref{fast}) and, for $t=u+iv$,
let $G(e^{-t})=a\log{\frac{1}{t}} +O(f(\mid t\mid))$ as $t\to 0$,
where $u>0$, $v=O(u^{1+\epsilon})$ as $u\to 0^+$, where $\epsilon$
and $a$ are as in part (i). Then, we have
$$
\frac{1}{p(n)} x^n[P(x)G(x)]
=a\log{\left(\frac{\sqrt{24n-1}}{2\pi}\right)} +O(n^{-1/2}
+f(c/\sqrt{n}+O(n^{-1/2-\epsilon}))
$$
with $c$ given by (\ref{c}).
\end{lemma}

As in \cite{GKW14}, we remark that parts (i) and (ii) can be combined
so that Lemma 2 generalizes to mixed asymptotic expansions
involving sums of powers of $t$ and logarithms of $1/t$. The proof
of Lemma 2, based on the saddle point method, is presented in
\cite[Section~3]{GKW14}.

\section{Proofs}

{\it Proof of Theorem 1.}

First, we set $m=1$ in Lemma 1(i). We have
\begin{equation}\label{misone}
\sum_{n=1}^\infty \mathbb{P}(M_n=1)x^n =xP(x).
\end{equation}
This implies at once (\ref{mnexact}). The asymptotic behavior of
the quotient in (\ref{mnexact}) may be found using Rademacher's
"exact-asymptotic" formula \cite{R37} (see also \cite[Chapter~5]{A76}). It
seems that a quicker way is to apply the result of Lemma 2(i). Here we have
$G(x)=x$, which obviously satisfies (\ref{fast}). Setting
$x=e^{-t}$ and expanding $e^{-t}$ by Taylor formula, we can write there as many powers of $t$
as we wish, which will be transferred into powers of $n^{-1/2}$ in the
asymptotic expansion of $\mathbb{P}(M_n=1)$. We decide to bound the error of
estimation up to a term of order $O(n^{-3/2})$ and write
\begin{equation}\label{e}
e^{-t} =1-t+\frac{1}{2}t^2 +f(t),
\end{equation}
with
\begin{equation}\label{f}
f(t)=\sum_{j=3}^\infty\frac{t^j}{j!}.
\end{equation}
The representation (\ref{e}) requires to apply Lemma 2(i) twice:
for the term $-t$ with $a=-1$ and $b=1$ and for the term
$\frac{1}{2} t^2$ with $a=1/2$ and $b=2$. Furthermore, (\ref{f})
implies that $f(c/\sqrt{n}+O(n^{-1/2}-\epsilon))=O(n^{-3/2})$.
Thus, from (\ref{misone}) it follows that
\begin{equation}\label{mnexpand}
\mathbb{P}(M_n=1) =\frac{x^n[xP(x)]}{p(n)} =1-A_1(n)+A_2(n)
+O(n^{-3/2}).
\end{equation}
The computation of $A_1(n)$ and $A_2(n)$ is based on Lemma 2(i)
with $s$ given by (\ref{es}). We have
\begin{eqnarray}\label{first}
& & A_1(n) =\frac{2\pi}{\sqrt{24n-1}}\frac{s}{s-1} \sum_{j=0}^2
\frac{(2+j)!}{j!(2-j)!}\left(-\frac{1}{2s}\right)^j \nonumber \\
& & =\frac{c}{\sqrt{n}}\left(1-\frac{1}{24n}\right)^{-1/2}
\left(1-\frac{2}{s}+\frac{1}{s^2}+O(s^{-3}\right) O(e^{-2s})
\nonumber \\
& & =\frac{c}{\sqrt{n}}\left(1+\frac{1}{48n}+O(n^{-2})\right)
\left(1-\frac{1}{c\sqrt{n}}\left(1-\frac{1}{24n}\right)^{-1/2}
+O(n^{-1})\right) \\
& & =-\frac{c}{\sqrt{n}}(1+O(n^{-1}))\left(1-\frac{1}{c\sqrt{n}}
+O(n^{-1})\right) =-\frac{c}{\sqrt{n}}-\frac{1}{n}+O(n^{-3/2}).
\nonumber
\end{eqnarray}
In the same way, for $A_2(n)$, we obtain
\begin{equation}\label{second}
A_2(n) =\frac{c^2}{2n} +O(n^{-3/2}).
\end{equation}
The proof is completed by substituting (\ref{first}) and
(\ref{second}) into (\ref{mnexpand}). $\Box$

{\it Proof of the Corollary.}

The total probability formula and the asymptotic estimate given in Theorem 1
imply that
\begin{eqnarray}
& & \mathbb{P}\left(\frac{L_n M_n}{\sqrt{n}}-\frac{1}{2c}
\log{n}\le u \right) =\mathbb{P}\left(\frac{L_n
M_n}{\sqrt{n}}-\frac{1}{2c} \log{n}\le u \mid M_n=1\right)
\mathbb{P}(M_n=1) \nonumber \\
& & +\mathbb{P}\left(\frac{L_n M_n}{\sqrt{n}}-\frac{1}{2c}
\log{n}\le u\mid M_n\neq 1\right) \mathbb{P}(M_n\neq 1) \nonumber \\
& & =\mathbb{P}\left(\frac{L_n}{\sqrt{n}} -\frac{1}{2}\log{n} \le u\right)(1+O(1/\sqrt{n}))
+O(1/\sqrt{n}). \nonumber
\end{eqnarray}
Hence the Corollary follows easily from Erd\"{o}s and Lehner's result (\ref{el}). $\Box$

{\it Proof of Theorem 2.}

First, we represent the function $F(x)$ given by (\ref{efx}) as
\begin{equation}\label{efonetwo}
F(x)=F_1(x)+F_2(x),
\end{equation}
where
\begin{equation}\label{efone}
F_1(x) =\sum_{k=1}^\infty kx^k\prod_{j=k+1}^\infty(1-x^j),
\end{equation}

\begin{eqnarray}\label{eftwo}
& & F_2(x) =\sum_{k=1}^\infty k(\sum_{l=2}^\infty x^{kl})
\prod_{j=k+1}^\infty (1-x^j) \nonumber \\
& & =\sum_{k=1}^\infty \frac{kx^{2k}}{1-x^k} \prod_{j=k+1}^\infty
(1-x^j).
\end{eqnarray}
Grabner and Knopfmacher \cite[formula~6.2]{GK06} found a simple
alternative representation for $F_1(x)$:
$$
F_1(x)=\sum_{k=1}^\infty\frac{x^k}{1-x^k}.
$$
It is also known that
$$
\sum_{n=1}^\infty p(n)\mathbb{E}(L_n)x^n =P(x)F_1(x)
$$
(see, e.g., \cite[p.~1059]{GKW14}). In addition, Grabner et al. \cite[p.~1084]{GKW14} used Mellin trasform technique to show that
$$
F_1(e^{-t}) =\frac{\log{(1/t)}+\gamma}{t}+\frac{1}{4}
-\frac{t}{144} +O(\mid t\mid^3), \quad t\to 0.
$$
From this expansion and their main result (see also both parts of
Lemma 2) they derived asymptotic formula (\ref{eel}) for
$\mathbb{E}(L_n)$. From (\ref{efonetwo}) it follows that
$x^n[F(x)]=x^n[F_1(x)]+x^n[F_2(x)]$, which in turn implies that
\begin{equation}\label{elm}
\mathbb{E}(L_n M_n) =\mathbb{E}(L_n) +\frac{x^n[F_2(x)]}{p(n)}.
\end{equation}

The asymptotic analysis of the second summand in the right-hand
side of (\ref{elm}) will be based on Lemma 2(ii). So, we need a
suitable expansion for $F_2(e^{-t})$. We set $t=u+iv$ with $u$ and
$v$ satisfying the conditions of Lemma 2(ii) and focus on an
asymptotic estimate for $F_2(e^{-u})$ as $u\to 0^+$. For the sake
of convenience, we also set
\begin{equation}\label{gk}
g_k(u) =\prod_{j=k+1}^\infty (1-e^{-ju}), \quad u>0.
\end{equation}
Representing $g_k(u)$ as a Riemann sum with step size $u\to 0^+$
and replacing it by the corresponding integral, we obtain the
following basic estimate:
\begin{equation}\label{gkbas}
g_k(u)
=\exp{(\sum_{j>k}\log{(1-e^{-ju})})}=\exp{\left(\frac{1}{u}\int_{ku}^\infty
\log{(1-e^{-w})}dw+o(1)\right)}.
\end{equation}
Next, we proceed  with the representation:
\begin{equation}\label{sum}
F_2(e^{-u})=\sum_{j=1}^4 S_j(u),
\end{equation}
where by (\ref{eftwo}) and (\ref{gk})
\begin{equation}\label{sj}
S_j(u) =\sum_{k\in I_j}\frac{ke^{-2ku}}{1-e^{-ku}} g_k(u), \quad
j=1,2,3,4.
\end{equation}
The subintervals $I_j$ are defined as follows. Let $\alpha>0$ be
fixed. We set
\begin{equation}\label{ione}
I_1=[1,\alpha/u],
\end{equation}
\begin{equation}\label{itwo}
I_2=\left(\frac{\alpha}{u},
\frac{1}{u}\left(\log{\frac{1}{u}}-\log{\log{\left(\frac{1}{u}\right)^3}}\right)\right],
\end{equation}
\begin{eqnarray}\label{ithree}
& & I_3 \nonumber \\
& &
=\left(\frac{1}{u}\left(\log{\frac{1}{u}}-\log{\log{\frac{1}{u}}}-\log{3}\right),
\frac{1}{u}\left(\log{\frac{1}{u}}+\log{\log{\frac{1}{u}}}+\log{2}\right)\right],
\end{eqnarray}
\begin{equation}\label{ifour}
 I_4
=\left(\frac{1}{u}\left(\log{\frac{1}{u}}
+\log{\log{\left(\frac{1}{u}\right)}}\right)+\log{2},\infty\right).
\end{equation}

We start with an estimate for $S_1(u)$. Since, for $k\in I_1$,
$$
\int_{ku}^\infty\log{(1-e^{-w})}dw \le\int_\alpha^\infty
\log{(1-e^{-w})}dw=-c_\alpha<0,
$$
from (\ref{gkbas}) we get
$$
g_k(u)=O(e^{-c_\alpha/u}).
$$
Hence, using (\ref{ione}), we obtain
\begin{eqnarray}\label{esone}
& & S_1(u) =O(e^{-c_\alpha/u})\sum_{u\le ku\le\alpha}
\frac{ke^{-2ku}}{1-e^{-ku}} \nonumber \\
& &
=O(u^{-2}e^{-c_\alpha/u})\int_0^\alpha\frac{we^{-2w}}{1-e^{-w}}dw
=O(u^{-2}e^{-c_\alpha/u}).
\end{eqnarray}

As a preparation for the estimate of $S_2(u)$, we first note that a
single integration by parts in the integral from the right-hand side of (\ref{gkbas}) yields
$$
\int_{ku}^\infty\log{(1-e^{-w})}dw =-ku\log{(1-e^{-ku})}
-\int_{ku}^\infty\frac{w}{e^w-1}dw.
$$
On the other hand, from formula 27.1.2 in \cite{AS65} we have
\begin{eqnarray}
& & \int_y^\infty\frac{w}{e^w-1}dw =y\sum_{k=1}^\infty\frac{e^{-ky}}{k}
+\sum_{k=1}^\infty\frac{e^{-ky}}{k^2} \nonumber \\
& & =-y\log{(1-e^{-y})} +e^{-y} +O(e^{-2y}), \quad y\to\infty. \nonumber
\end{eqnarray}
Applying this estimate to (\ref{gkbas}), we obtain
\begin{equation}\label{gksec}
g_k(u) =\exp{\left(-\frac{1}{u}e^{-ku} +O(u^{-1}e^{-2ku}) +o(1)\right)}
\end{equation}
for all $k$ that satisfy $ku\to\infty$ as $u\to 0^+$. So, from (\ref{itwo}) it follows that
\begin{eqnarray}\label{estwo}
& & S_2(u) =\frac{1+o(1)}{u} \sum_{k\in I_2}\frac{ku e^{-2ku}}{1-e^{-ku}}
\exp{\left(-\frac{1}{u}e^{-ku} +O(u^{-2}e^{-2ku})\right)} \nonumber \\
& & =O\left(\frac{\log{\frac{1}{u}}}{u}\left(\sum_{k\in I_2}\frac{e^{-2\alpha}}{1-e^{-\alpha}}\right)
\exp{\left(-\log{\left(\frac{1}{u}\right)^3} +O\left(\frac{1}{u}e^{-2\log{\frac{1}{u}}}\right)\right)}\right) \nonumber \\
& & =O\left(\frac{u^3\log^2{\frac{1}{u}}}{u^2}\right) =O\left(u\log^2{\frac{1}{u}}\right),\quad u\to 0^+.
\end{eqnarray}
(In the third line of (\ref{estwo}) we have also used that the sum over $I_2$ contains
at most $\frac{1}{u}\log{\frac{1}{u}}$ summands.)

We proceed now to the estimate for $S_3(u)$, whose contribution to $F_2(e^{-u})$ is
the most essential one. First, it is easy to
observe that
$$
0\le\frac{1}{1-e^{-ku}} =O\left(u\log{\frac{1}{u}}\right), \quad u\to 0^+
$$
uniformly for all $k\in I_3$. Hence, approximating once again a Riemann sum
by the corresponding integral, we obtain
\begin{eqnarray}\label{esthreeasone}
& & S_3(u) =\left(1+O\left(u\log{\frac{1}{u}}\right)\right) u^{-2}
\int_{\log{\frac{1}{u}}-\log{\log{\frac{1}{u}}}-\log{3}}^{\log{\frac{1}{u}}+\log{\log{\frac{1}{u}}}+\log{2}}
we^{-2w} e^{-\frac{1}{u}e^{-w}} dw \nonumber \\
& & =\left(1+O\left(u\log{\frac{1}{u}}\right)\right)
\int_{-\log{\log{\frac{1}{u}}}-\log{3}}^{\log{\log{\frac{1}{u}}}+log{2}}
\left(z+\log{\frac{1}{u}}\right) e^{-2z} e^{-e^{-z}} dz.
\end{eqnarray}
Changing the variable in ({\ref{esthreeasone}), we represent $S_3(u)$ as follows:
\begin{equation}\label{esthreeastwo}
S_3(u)=(J_1(u)+J_2(u))\left(1+O\left(u\log{\frac{1}{u}}\right)\right),
\end{equation}
where
\begin{equation}\label{jone}
J_1(u)=\int_{2/\log{\frac{1}{u}}}^{3\log{\frac{1}{u}}} (-\log{y})ye^{-y}dy,
\end{equation}

\begin{equation}\label{jtwo}
J_2(u) =\left(\log{\frac{1}{u}}\right) \int_{2/\log{\frac{1}{u}}}^{3\log{\frac{1}{u}}}
ye^{-y} dy.
\end{equation}
An easy calculation based on an asymptotic estimate for the incomplete gamma function
\cite[formula~6.5.32]{AS65} shows that
\begin{eqnarray}
& & \int_{2/\log{\frac{1}{u}}}^{3\log{\frac{1}{u}}}
ye^{-y} dy =\int_0^\infty ye^{-y} dy -\int_0^{2/\log{\frac{1}{u}}} ye^{-y} dy
-\int_{3\log{\frac{1}{u}}}^\infty ye^{-y}dy \nonumber \\
& & =1-O\left(\left(\log{\frac{1}{u}}\right)^{-2}\right)
-O\left(\left(\log{\frac{1}{u}}\right)e^{-3\log{\frac{1}{u}}}\right)
= 1- O\left(\left(\log{\frac{1}{u}}\right)^{-2}\right). \nonumber
\end{eqnarray}
Inserting this estimate into (\ref{jtwo}), we obtain
\begin{equation}\label{jtwoas}
J_2(u)=\log{\frac{1}{u}} +O\left(1/\log{\frac{1}{u}}\right).
\end{equation}
For the estimate of $J_1(u)$, we recall that
$$
\int_0^\infty ye^{-y}(\log{y}) dy=1-\gamma
$$
(see, e.g., \cite[formula~865.902]{D61}. Hence we have
$$
J_1(u) =\gamma-1 -\int_0^{2/\log{\frac{1}{u}}} ye^{-y}(-\log{y})dy
-\int_{3\log{\frac{1}{u}}}^\infty ye^{-y}(-\log{y})dy.
$$
Both integrals in the right-hand side of this representation are negligible.
The first one can be estimated bounding $e^{-y}$ by $1$ and then integrating by parts.
This gives a bound of order $O\left(\log{\log{\frac{1}{u}}}/\log^2{\frac{1}{u}}\right)$.
For the second one we can use once again \cite[formula~6.5.32]{AS65},
which leads to the upper bound $O\left(u^3\log^2{\frac{1}{u}}\right)$. Thus we
have
\begin{equation}\label{joneas}
J_1(u) =\gamma-1 +O\left(\left(\log{\log{\frac{1}{u}}}\right)/\log^2{\frac{1}{u}}\right).
\end{equation}
Combining (\ref{esthreeastwo}) - (\ref{joneas}), we finally obtain
\begin{equation}\label{esthreeasthree}
S_3(u)=\log{\frac{1}{u}}+\gamma-1 +O\left(\frac{1}{\log{\frac{1}{u}}}\right).
\end{equation}

We end up our analysis with an estimate for $S_4(u)$. First, we bound the last
exponent in the right-hand side of (\ref{gksec}) by $1$ and then, as previously,
we approximate $S_4(u)$ by an integral, which can be bounded using again \cite[formula~6.5.32]{AS65}. Then, in a similar way, from (\ref{ifour}) we obtain
\begin{eqnarray}\label{esfouras}
& & S_4(u) =O\left(u^{-2}\int_{\log{\frac{1}{u}}+\log{\log{\frac{1}{u}}}+\log{2}}^\infty
we^{-2w} dw\right) \nonumber \\
& &= O\left(\int_{\log{\log{\frac{1}{u}}}+\log{2}}^\infty \left(z+\log{\frac{1}{u}}\right)
e^{-2z} dz\right) \nonumber \\
& & =O\left(\left(\log{\log{\frac{1}{u}}}\right)/\log^2{\frac{1}{u}}\right) +
O\left(1/\log{\frac{1}{u}}\right) =O\left(1/\log{\frac{1}{u}}\right).
\end{eqnarray}
Now, from (\ref{sum}), (\ref{esone}), (\ref{estwo}), (\ref{esthreeasthree}) and
(\ref{esfouras}) it follows that
\begin{eqnarray}\label{eftwoas}
& & F_2(e^{-u}) =O(u^{-2}e^{-c_\alpha/u}) +O\left(u^2\log{\frac{1}{u}}\right)
+\log{\frac{1}{u}}+\gamma-1 +O\left(1/\log{\frac{1}{u}}\right) \nonumber \\
& & =\log{\frac{1}{u}}+\gamma-1 +O\left(1/\log{\frac{1}{u}}\right).
\end{eqnarray}

To transfer the variable $u$ in (\ref{eftwoas}) into $t=u+iv$, we recall the
relationship between $u$ and $v$ in Lemma 2(ii). First, we note that
\begin{equation}\label{logtu}
\log{\frac{1}{t}}=\log{\frac{1}{u}} +\log{\frac{1}{1+iv/u}} =\log{\frac{1}{u}}
+O(\mid v\mid/u) =\log{\frac{1}{u}} +O(u^\epsilon)
\end{equation}
and
\begin{equation}\label{etu}
e^{-t}=e^{-u}e^{-iv} =e^{-u}(1+O(\mid v\mid)=e^{-u}+O(u^{1+\epsilon})
\end{equation}
as $u\to 0^+$.
Hence, by Taylor formula,
\begin{equation}\label{taylor}
F_2(e^{-t}) =F_2(e^{-u}) +O\left(u^{1+\epsilon}\mid\frac{d}{du}F_2(e^{-u})\mid\right).
\end{equation}
The last error term can be estimated following the same line of reasoning.
We outline briefly the proof. By (\ref{eftwo}) and (\ref{gk})
we have
\begin{eqnarray}\label{eftwoprime}
& & \frac{d}{du}F_2(e^{-u}) =-\sum_{k=1}^\infty\frac{k^2 e^{-ku}}{1-e^{-ku}}g_k(u)
-\sum_{k=1}^\infty\frac{k^2 e^{-ku}}{(1-e^{-ku})^2}g_k(u) \nonumber \\
& & -\sum_{k=1}^\infty\frac{k e^{-ku}}{1-e^{-ku}}g_k^\prime(u),
\end{eqnarray}
where
$$
g_k^\prime(u)=g_k(u)\sum_{j=k+1}^\infty\frac{je^{-ju}}{1-e^{-ju}}.
$$
The main contribution to the asymptotic of $\mid\frac{d}{du}F_2(e^{-u})\mid$
as $u\to 0^+$ is given by the first term of the right-hand side of (\ref{eftwoprime})
(since it contains the factor $k^2$). If we break up the range of summation as
previously into the union $\bigcup_{j=1}^4 I_j$ (see (\ref{ione}) - ({\ref{ifour})), we
can again conclude that the contribution of the sum over $I_3$ is the largest one.
An approximation by a Riemann integral requires division and multiplication by the
cube of the step size, $u^3$. Hence, instead of $u^{-2}$, the factor $u^{-3}$ will
multiply the same integral from the right-hand side of (\ref{esthreeasone}).
An argument similar to (\ref{esthreeasthree}) and (\ref{eftwoas}) implies
the following bound:
$$
\frac{d}{du}F_2(e^{-u}) =O\left(\frac{1}{u}\log{\frac{1}{u}}\right), \quad u\to 0^+.
$$
Consequently, the remainder term in (\ref{taylor}) is $O\left(u^\epsilon\log{\frac{1}{u}}\right)$.
Combining this with (\ref{eftwoas}) - (\ref{eftwoprime}), we obtain
$$
F_3(e^{-t}) =\log{\frac{1}{t}}+\gamma -1+ O\left(1/\log{\frac{1}{t}}\right), \quad t\to 0.
$$
Therefore, we are ready to apply Lemma 2(ii) with $G(x):=F_2(x)$ and $f(t)=1/\log{\frac{1}{t}}$.
We obtain
\begin{equation}\label{as}
\frac{1}{p(n)} x^n[F_2(x)] =\log{\frac{c}{\sqrt{n}}} +\gamma-1 +O(1/\log{n}), \quad n\to\infty.
\end{equation}
Now, we recall (\ref{elm}), which in combination with (\ref{as}) completes the proof of Theorem 2. $\Box$

\section{Concluding Remarks}

The main goal of this study is the comparison between the typical growths of the first block area $L_n M_n$ and
its base length $L_n$ in the Ferrers diagram of a random integer partition $n$. It turns out that the leading terms
in the asymptotic expansions of the expectations of these two statistics are the same for large $n$; both are equal to
$\frac{\sqrt{n}}{2c}\log{n}$. Erd\"{o}s and Lehner's limit
theorem (\ref{el}) and Theorem 1 show that this leading term
controls the weak convergence of $L_n$ and $L_n M_n$ both tending to a Gumbel distributed random variable
after one and the same appropriate normalization.
The expectations of $L_n$ and $L_n M_n$ are, however, different for large $n$. In fact, by Theorem 2,
that
$$
\lim_{n\to\infty}\left(\mathbb{E}(L_n M_n)-\mathbb{E}(L_n)-\frac{1}{2}\log{n}\right) =-C =-0.67165... .
$$
This phenomenon suggests a question related to the shape of a random Ferrers diagram of $n$.
It is interesting to determine how far from the largest block area the $L_n M_n$-area is if all block areas of the Ferrers diagram
are arranged in non-increasing order. This kind of rearrangement of blocks of the Ferrers diagram was studied by Fristedt \cite[p.~710]{F93},
who obtained a general limit theorem for the $r$th block area. To state the problem in more precise way, we
denote by $X_n^{(k)}$ the multiplicity of part $k$ ($k=1,...,n$) in a random integer partition of $n$. Let $Z_n^{(r)}$ be
the $r$th largest member of the sequence $\{kX_n^{(k)}\}_{k=1}^n$. Erd\"{o}s and Szalay \cite{ES84} showed that
$$
\lim_{n\to\infty}\mathbb{P}\left(\frac{c}{\sqrt{n}}Z_n^{(1)} -\frac{1}{2}\log{\frac{n}{c^2}} -\log{\log{\log{n}}}\le u\right)
=e^{-e^{-u}}, \quad -\infty<u<\infty.
$$
Fristedt \cite[Theorem~2.7]{F93} generalized this result to $Z_n^{(r)}$, where $r\ge 1$ is fixed. He showed that
\begin{eqnarray}\label{fri}
& & \lim_{n\to\infty}\mathbb{P}\left(\frac{c}{\sqrt{n}}Z_n^{(r)} -\frac{1}{2}\log{\frac{n}{c^2}} -\log{\log{\log{n}}} \le u\right)
\nonumber \\
& & =\int_{-\infty}^u \frac{\exp{(-e^{-w}-r w)}}{(r-1)!} dw, \quad -\infty<u<\infty.
\end{eqnarray}
Let $R_n$ be the smallest value of $r$ such that $Z_n^{(r)}=L_n M_n$. We conjecture that
\begin{equation}\label{con}
\mathbb{E}(R_n)\asymp\log{\log{n}}, \quad n\to\infty.
\end{equation}
This claim is supported by the following heuristic but non-rigorous argument.
From the result of Theorem 2 it follows that
\begin{equation}\label{t}
\mathbb{E}(L_n M_n) =\frac{\sqrt{n}}{2c}\log{n} +O(\sqrt{n}), \quad n\to\infty.
\end{equation}
In addition, a calculation of the expectation of the distribution
in the right-hand side of (\ref{fri}) given by Fristedt \cite[p.~708]{F93} could be used to show that
\begin{equation}\label{ez}
\mathbb{E}(Z_n^{(r)}) =\frac{\sqrt{n}}{2c}\left(\log{n} +2\log{\log{\log{n}}} -2\log{r}\right) +O(\sqrt{n}),
\end{equation}
if $r=r(n)\to\infty$ as $n\to\infty$. Combining (\ref{t}) with (\ref{ez}), one may conclude
that $\log{r(n)}\approx \log{\log{\log{n}}}$,
which supports the claim in (\ref{con}). We hope to return to this question in a future study.


\begin{thebibliography}{15}

\bibitem{AS65}
Abramovitz, M., Stegun, I.A.: Handbook of Mathematical Functions
with Formulas, Graphs and Mathematical Tables. Dover Publ., Inc.
New York (1965)

\bibitem{A76}
Andrews, G.E.: The Theory of Partitions. Encyclopedia Math. Appl.
2, Addison-Wesley, Reading, MA (1976)

\bibitem{ABBKM16}
Archibald, M., Blecher, A., Brennan, C., Knopfmacher, A., Mansour,
T.: Partitions according to multiplicities and part sizes. Australas.
J. Combin. {\bf 66}, 104-119 (2016)

\bibitem{D61}
Dwight, H.B.: Tables of Integrals and Other Mathematical Data, 4th ed.
MacMillan, New York (1961)

\bibitem{EL41}
Erd\"{o}s, P., Lehner, J.: The distribution of the number of summands
in the partitions of a positive integer. Duke Math. J. {\bf 8}, 335-345
(1941)

\bibitem{ES84}
Erd\"{o}s, P., Szalay, M.: On the statistical theory of partitions. In:
Topics in Classical Number Theory, vol. I (G. Halasz ed.). North-Holland,
Amsterdam, pp.397-450 (1984)

\bibitem{F93}
Fristedt, B.: The structure of random partitions of large integers. Trans.
Amer. Math. Soc. {\bf 337}, 703-735 (1993)

\bibitem{GK06}
Grabner, P., Knopfmacher, A.: Analysis of some new partition statistics.
Ramanujan J. {\bf 12}, 439-454 (2006)

\bibitem{GKW14}
Grabner, P., Knopfmacher, A., Wagner, S.: A general asymptotic
scheme for the analysis of partition statistics. Combin. Probab.
Comput. {\bf 23}, 1057-1086 (2014)

\bibitem{HR18}
Hardy, G.H., Ramanujan S D.: Asymptotic formulae in combinatory
analysis. Proc. London Math. Soc. {\bf 17(2)}, 75-115 (1918)

\bibitem{H38}
Husimi, K.: Partitio numerium as occurring in a problem of
nuclear physics. Proc. Phys.-Math. Soc. Japan {\bf 20}, 912-925 (1938)

\bibitem{KL76}
Kessler, I., Livingston, M.: The expected number of parts in a partition
of $n$. Monatsh. Math. {\bf 81}, 203-212 (1976)

\bibitem{R37}
Rademacher, H.: On the partition function $p(n)$. Proc. London Math, Soc.
{\bf 43}, 241-254 (1937)

\bibitem{R74}
Richmond, L.B.: The moments of partitions, I. Acta Arith. {\bf 211}, 345-373
(1912)

\bibitem{S}
Sloane, N.: On-line Encyclopedia of Integer Sequences (https://oeis.org/)
\end{thebibliography}
\end{document}